\documentclass[11pt,a4paper]{article}

\usepackage{amsmath}
\usepackage{amsthm}
\usepackage{amsfonts}
\usepackage[koi8-r]{inputenc}
\usepackage[russian]{babel}

\parskip\smallskipamount
\newtheorem{theorem}{Theorem}
\theoremstyle{remark}
\newtheorem{property}{Property}
\newtheorem{rem}{Remark}

\renewcommand{\Pr}{\mathbf{P}}
\newcommand{\E}{\mathbf{E}}

\newcommand{\ov}[1]{\overline{#1}}
\newcommand{\LL}[1]{\mathcal{L}_#1}
\renewcommand{\SS}[1]{\mathcal{S}_#1}
\renewcommand{\phi}{\varphi}
\renewcommand{\epsilon}{\varepsilon}

\parskip \smallskipamount

\title{On the exact asymptotics for the stationary sojourn time distribution
in a tandem of queues with light-tailed service times 
\footnote{The research was partially supported by the 
EPSRC Grant EP/E033717/1}}

\author{S.G.~Foss\\
\emph{Heriot-Watt University, Edinburgh and}\\
\emph{Institute of Mathematics, Novosibirsk}}
\date{}

\begin{document}

\maketitle

\begin{quotation}\small
We study the asymptotics of the stationary sojourn time $Z$ of a ``typical 
customer'' in a tandem of single-server queues. It is shown that,
in a certain ``intermediate'' region of light-tailed service time
distributions, $Z$ may take a large value mostly due to a 
large value of a single service time of one of customers. 
Arguments used in the paper allow us to obtain also an elementary proof of
the logarithmic asymptotics for the tail distribution of the stationary sojourn time
in the whole class of light-tailed distributions. 
\end{quotation}

{\bf Keywords:} tandem of queues, sojourn time, large deviations,
Cramer condition, exact and logarithmic asymptotics, class 
$\SS{\gamma}$.

\section{Introduction and main results}
\label{sec:introduction}
Consider an open queueing network which is a tandem of two single-server queues 
$GI/GI/1 \rightarrow /GI/1$ with ``first-come-first-served'' service disciplines.

Consider three mutually independent sequences of non-negative random
variables
$\{ \tau_n\}$, $\{\sigma_n^{(1)} \}$ 
and $\{\sigma_n^{(2)} \}$, each of which is an i.i.d. sequence.
Here $\tau_n$ is the inter-arrival time between the $(n-1)$st and $n$th
customers 
(with mean $a= {\mathbf E} \tau_1$). 
Customer $n$ receives service in the first queue
of duration $\sigma_n^{(1)}$ ( with 
distribution function $G^{(1)}$ and positive mean
$b^{(1)} ={\mathbf E} \sigma_1^{(1)}$) and then in
the second of duration  
  $\sigma_n^{(2)}$ (with distribution function $G^{(2)}$ and
positive mean $b^{(2)} ={\mathbf E} \sigma_1^{(2)}$).
Is it assumed that the network is stable, i.e.
$\max \left(b^{(1)}, b^{(2)}\right) < a$.
It is well-known (see, for example, \cite{Lo}) that, under this
assumption, there exists a unique stationary (limiting)
distribution of the sojourn time $Z$ in the network (i.e. of
the duration of time between a customer's arrival to the first queue
and its departure from the second queue) and, for any initial condition,
the distribution of the sojourn time $Z_n$ of customer $n$ converges
in the total variation norm to the limiting distribution, 
as $n\to\infty$.

The following representation is also known (see, e.g., \cite{BF94}):
\begin{equation}\label{main}
Z = \sup_{0\le n \le m < \infty}
\left(
\sum_{-m}^{-n} \sigma_j^{(1)} + \sum_{-n}^0 \sigma_j^{(2)} -
\sum_{-m}^{-1} \tau_j
\right).
\end{equation}
One may interpret formula (\ref{main}) as follows.
Assume that the network has been working for an infinitely long time, 
starting
from time  $-\infty$. Then
$Z=Z_0$ is the sojourn time of customer $0$ that arrives at time instant $t=0$. 
Note that 
 $Z$ is a monotone function of all variables in the right-hand side
of \eqref{main}: it monotonically increases with the growth of any of the $\sigma$'s 
and with the decrease of any of the $\tau$'s.

The paper deals with a study of the asymptotics of probability
${\mathbf P}(Z>x)$  
when $x$ tends to infinity.
We consider the case where service times distributions in both queues
are light-tailed, i.e.
\begin{equation}\label{cond1}
 \phi_{\sigma_1^{(i)}} (\lambda ) <\infty
\end{equation}
for $i=1,2$ and for some positive $\lambda$. Here we use the
standard notation: 
$\phi_X (\lambda) = {\mathbf E} e^{\lambda X}$ 
is the exponential moment of a random variable 
$X$ at point $\lambda$. For short, we write in
the sequel 
$\phi_{(i)}(\lambda) = \phi_{\sigma_1^{(i)}}(\lambda)$ for
$i=1,2$ and $\phi_{\tau}(\lambda) = \phi_{\tau_1}(\lambda)$.
For $i=1,2$, denote 
\begin{displaymath}
\gamma^{(i)} = \sup \{\lambda \ : \ \phi_{(i)}
(\lambda) \phi_{\tau}(-\lambda) \le 1 \}\in (0,\infty )
\end{displaymath}
and let 
$$\gamma = \min \left( \gamma^{(1)}, \gamma^{(2)} \right).$$
For two positive functions
$f_1$ and $f_2$ and for a constant $d\ge 0$, the notation 
$f_1(x) \sim df_2(x)$ means that
$f_1(x)/f_2(x) \to d$ as $x\to\infty$. In particular, if
$d=0$, then $f_1(x) = o(f_2(x))$. 

The following logarithmic (``rough'') asymptotics holds 
(see, for example, \cite{BPT, GAN}):  
\begin{theorem} \label{th1}
Under condition \eqref{cond1}, 
$$
- \ln {\mathbf P} (Z>x) \sim \gamma x.
$$
\end{theorem}
All (known to us) proofs of Theorem \ref{th1} (see, for example,
\cite{BPT, GAN}) use the techniques of large deviations.

Our main result (Theorem 2) provides the {\it exact}
asymptotics of large deviations, under the following additional
assumption:
\begin{equation} \label{eq:11}
R \equiv \max \left(
 \phi_{(1)}
(\gamma) , \phi_{(2)}
(\gamma)  
\right) \phi_{\tau}(-\gamma) < 1.
\end{equation}
In particular, (\ref{eq:11}) implies the finiteness of ${\mathbf E} e^{\gamma Z}$.
Indeed, 
\begin{equation}\label{konZ}
{\mathbf E} e^{\gamma Z} \le \sum_{0\le n\le m }
{\mathbf E} \exp \left( \gamma \left(
\sum_{-m}^{-n} \sigma_j^{(1)} + \sum_{-n}^0 \sigma_j^{(2)} -
\sum_{-m}^{-1} \tau_j \right)
\right) \le (1-R)^{-2} \varphi_{\tau}^{-2}(-\gamma) < \infty.
\end{equation}
In order to state Theorem 2 , we need a number of further definitions
and notation. 

We use the same symbol
$F$ to denote a probability distribution on
the real line and also its distribution function. 
Let $\ov{F}$ be the tail of distribution
$F$, i.e..
$\overline{F}(x)=1-F(x)$, and $F^{*n}$ the 
$n$-fold convolution of $F$.
A distribution function~$F$
{\it belongs to the class}~$\LL{\beta}$, $\beta\ge0$ if  
\begin{equation}
  \label{eq:2}
    \text{$\ov{F}(x)>0$ for all $x$} \quad\text{and} 
    \lim_{x\to\infty}\frac{\ov{F}(x-y)}{\ov{F}(x)} = e^{\beta y} 
    \quad\text{for any fixed $y$}.
\end{equation}
Due to the monotonicity of $\overline{F}$, 
the convergence in \eqref{eq:2} is necessarily uniform in
$y$ on any compact set. Therefore we may find such a function
$h(x) \uparrow \infty$, $h(x) = o(x)$ that
\begin{equation} \label{hx}
\lim_{x\to\infty} \sup_{|y| \le h(x)}
\left|
\frac{\overline{F}(x+y)}{\overline{F}(x)} e^{\beta y} - 1
\right| = 0.
\end{equation}
If $h_1$ and $h_2$ are two functions satisfying (6), 
then the function $h_3(x) = h_1(x)+h_2(x-h_1(x))$ has the same property.

The distribution function ~$F$ of a random variable $X$
{\it belongs to the class}~$\SS{\beta}$, 
$\beta\ge0$, if $F\in\LL{\beta}$,
$\phi_X(\beta)<\infty$ and
\begin{equation}
  \label{eq:29}
  \ov{F^{*2}}(x) = {\bf P} (X_1+X_2>x) \sim 2\phi_{X}(\beta)\ov{F}(x) 
  \qquad\text{as $x\to\infty$}
\end{equation}
where $X_1$ and $X_2$ are two independent copies of $X$.
Here, with necessity, 
\begin{eqnarray*}
{\bf P}(X_1+X_2>x) &\sim & 
{\bf P}(X_1+X_2>x, X_1\le h(x)) +
 {\bf P}(X_1+X_2>x, X_2\le h(x))\\
&\sim &
{\bf P}(X_1+X_2>x, X_1\ge x-h(x)) +
 {\bf P}(X_1+X_2>x, X_2\ge x-h(x))
\end{eqnarray*}
where $h(x)$ is any function satisfying condition \eqref{hx} (see,
for example, \cite{FZ}).

A typical example of a distribution $F\in \SS{\beta}$ with $\beta >0$ is
a distribution with the tail $\ov{F}(x) =
Cx^{-\alpha}e^{-\beta x}$  for some $\alpha >1$, $C\in (0,1]$
and all $x\ge 1$.

\begin{theorem} \label{th2}
Assume that condition \eqref{eq:11} holds. Suppose that
\begin{equation}\label{th2-f1}
 \ov{G^{(i)}}(x) \sim c^{(i)} \ov{F}(x), \quad i=1,2
\end{equation}
for some function $F\in \SS{\gamma}$
and constants $c^{(1)} \ge 0$ and $c^{(2)} \ge 0$. Then 
\begin{equation}\label{th2-f2}
{\mathbf P}(Z>x) \sim
\sum_{i=1}^2 \sum_{j=-\infty}^0
\Pr (Z>x, \sigma_j^{(i)} > x-h(x)) \sim K \ov{F}(x)
\end{equation}
where $h(x)$ is any function satisfying conditions
\eqref{hx} and the constant $K$ is determined below in formula 
(\ref{zvezda}).  
\end{theorem}
\begin{rem}
Since the exact representation of constant $K$ is complex and depends on characteristics
which are not ``computable'', it seems to be reasonable to provide useful
upper and lower bounds for $K$. For that, let, for $i=1,2$,
$R_i = \varphi_{(i)}(\gamma) \varphi_{\tau}(-\gamma)$
and
$R=\max (R_1,R_2)$. 
Then
\begin{equation}\label{up-low}
\frac{c^{(1)}\varphi_{(2)}(\gamma)}{1-R}
+\frac{c^{(2)}\varphi_{(1)}(\gamma)}{1-R_2} \le K \le
\frac{1}{(1-R_1)(1-R_2)}
\left( \frac{c^{(1)}\varphi_{(2)}(\gamma)}{1-R_1} + 
\frac{c^{(2)}\varphi_{(1)}(\gamma)}{1-R_2} \right).
\end{equation}
\end{rem}

\begin{rem}
The coefficients $c^{(1)}$ and $c^{(2)}$ in the statement of
Theorem \ref{th2} may be either positive or zero. If both coefficients are
positive, then, with necessity,
$\gamma^{(1)}=\gamma^{(2)}=\gamma$ and -- as it follows from Property 1 (see
Appendix) -- both distributions $G^{(i)}$,$i=1,2$, have to belong to the class
${\cal S}_{\gamma}$. 
If only one of the coefficients is positive, say if $c^{(1)}>0$ and
$c^{(2)}=0$, then the distribution $G^{(1)}$ belongs to the class 
${\cal S}_{\gamma}$ and $\gamma^{(2)}\ge \gamma^{(1)}=\gamma$. 
Finally, if $c^{(1)}=c^{(2)}=0$, then also $K=0$, as it follows from (\ref{up-low}).
\end{rem}

\begin{rem}
Taking into account the monotonicity properties (see the 
comments after formula (\ref{main}), we can obtain the following
``one-side'' analogues of Theorem \ref{th2}: \\
If, in the statement of Theorem \ref{th2}, one replaces condition (\ref{th2-f1}) by 
$\limsup_{x\to\infty} \overline{G^{(i)}}(x)/\overline{F}(x) \le c^{(i)}$, $i=1,2$
(or by
$\liminf_{x\to\infty} \overline{G^{(i)}}(x)/\overline{F}(x) \ge c^{(i)}$, $i=1,2$),
then the following holds:
$\limsup_{x\to\infty} \Pr (Z>x)/\overline{F}(x) \le K$ (or, respectively,
$\liminf_{x\to\infty} \Pr (Z>x)/\overline{F}(x) \ge K$), with the same 
$K$.
\end{rem}

\begin{rem}
A natural analogue of Theorem \ref{th2} 
holds for tandems of any finite number of queues and, more
generally, for tree-like queueing networks. However, an explicit representation
of the constant $K$ and even its bounds become less and less tractable as $n$
increases. Therefore we decided to consider the case $n=2$ only.
\end{rem}

\begin{rem}
The approach in the first part of the proof of Theorem
\ref{th2} (i.e. the construction of upper and lower bounds) allows us also
to obtain a simple proof of Theorem \ref{th1} without use of the techiques of the
large deviations theory (see Subsection 3.2). 
\end{rem}
\begin{rem}
The proposed method of proof is based on ideas developed in
 \cite{BF04} and applied therein to obtaining the distributional 
asymptotics for ${\mathbf P} (Z>x)$ in tandems of queues with subexponential
service times distributions (which are heavy-tailed).
Also, in \cite{BF04}, the asymptotics for the stationary waiting time in 
the second queue were obtained. Similar asymptotics may be found under the 
conditions of this paper, but by the use of essentially more complicated
formulae.
\end{rem}

\begin{rem}
In addition to the asymptotics for the tail $\Pr (Z>x)$, $x\to\infty$,
one can use results from \cite{FZ} to obtain the asymptotics 
for prestationary probabilities $\Pr (Z_k>x)$ (which are uniform in $k$). Here
$$
Z_k = \sup_{0\le n \le m \le k}
\left(
\sum_{-m}^{-n} \sigma_j^{(1)} + \sum_{-n}^0 \sigma_j^{(2)} -
\sum_{-m}^{-1} \tau_j
\right).
$$
\end{rem}

A proof of Theorem \ref{th2} is given in the next Section. Section 3 
contains a useful auxiliary information about the class ${\cal S}_{\beta}$,
a simple proof of Theorem \ref{th1} and some further comments.

\section{Proof}
\label{sec:proofs}

\subsection{Upper and lower bounds for the stationary sojourn time}
\label{sec:bounds}

{\bf Lower bound.} For $i=1,2$, let
\begin{equation}\label{Zi}
Z^{(i)} = \sigma_0^{(i)} + \max_{n\ge 0} 
\sum_{-n}^{-1} \left( \sigma_j^{(i)} - \tau_j \right)
\end{equation}
(where $\sum_0^{-1} = 0$). Then
$Z^{(1)} = \sigma_0^{(1)} + W^{(1)}$ where $W^{(1)} =
\max_{n\ge 0} 
\sum_{-n}^{-1} \left( \sigma_j^{(1)} - \tau_j \right)
$. Here  $Z^{(1)}$ (respectively, $W^{(1)}$) is the stationary 
sojourn (respectively, waiting) time in the first queue.
However, it is more convenient to propose a slightly different
interpretation of the formula above: 
$Z^{(1)}$ (respectively, 
$W^{(1)}$)
is the stationary sojourn (respectively, waiting) time
in an auxiliary tandem of queues where all service times in
the second queue are equal to zero.
Similarly, 
$Z^{(2)} = \sigma_0^{(2)} + W^{(2)}$ with $W^{(2)} =
\max_{n\ge 0} 
\sum_{-n}^{-1} \left( \sigma_j^{(2)} - \tau_j \right)
$,
and
$Z^{(2)}$ (respectively, $W^{(2)}$)
is the stationary sojourn (respectively, waiting) time in an
auxiliary tandem of queues where all service times in the first
queue and replaced by zero.

The monotonicity of $Z$  in all variables in \eqref{main} implies the
following bound
\begin{equation}\label{low1}
Z\ge \max \left( Z^{(1)}, Z^{(2)} \right) \quad\text{a.s.}
\end{equation}
and, in particular,
\begin{equation}\label{low3}
{\mathbf P} (Z>x) \ge \max ({\mathbf P} (Z^{(1)}>x), \ \  
{\mathbf P} (Z^{(2)}>x)).
\end{equation}

{\bf Upper bound.} Let $L\ge 1$ be an integer. Introduce an auxiliary single-server
queue with i.i.d. inter-arrival times $\widetilde{\tau}_n$ and (independent of them) i.i.d.
service times 
$\widetilde{\sigma}_n$ where
\begin{displaymath}
\widetilde{\tau}_n = \sum_1^L \tau_{(n-1)L+i}
\quad\text{and} \quad 
\widetilde{\sigma}_n = \max_{1\le j \le L}
\left(
\sum_{i=1}^j \sigma_{(n-1)L+i}^{(1)} +
\sum_{i=j}^L \sigma_{(n-1)L+i}^{(2)}
\right).
\end{displaymath}
Here the random variable $\widetilde{\sigma}_1$ may be viewed as
follows. Assume that customers $1,\ldots,L$ arrive simultaneously at
time instant $t=0$ into an empty network. Then $\widetilde{\sigma}_1$
is the time of the last completion of service of these customers in the
second queue.
It is not difficult to prove (see, for example, \cite{BF95}) that 
${\mathbf E}\widetilde{\sigma}_1/L
\to \max (b^{(1)}, b^{(2)})$ as $L\to\infty$. Hence   
${\mathbf E}\widetilde{\sigma}_1 < {\mathbf E} \widetilde{\tau}_1$
for all sufficiently large $L$. We fix such an  $L$ and define 
\begin{displaymath}
\widetilde{W} = \max_{n\ge 0} 
\sum_{-n}^{-1}
(\widetilde{\sigma}_n - \widetilde{\tau}_n) < \infty \quad\text{a.s.}
\end{displaymath}
(respectively $\widetilde{Z} = \widetilde{\sigma}_0+\widetilde{W}$),
the stationary sojourn (respectively waiting) time of customer $0$ in
this queueing system, which is a.s. finite. The monotonicity properties of
 $Z$ imply (see, for example, \cite{BF95, BF04}) that
\begin{equation}\label{up1}
Z \le \widetilde{Z} \quad\text{a.s.}
\end{equation}

\subsection{Proof of Theorem \ref{th2}.}

We prove Theorem \ref{th2} only for $c^{(1)}+c^{(2)}>0$, 
the statement in the case $c^{(1)}=c^{(2)}=0$ follows
by monotonicity.

>From [5, Theorem 1] and from inequality 
\eqref{low3}, we get 
\begin{equation}\label{low10}
\liminf_{x\to\infty} 
\frac{{\mathbf P} (Z>x)}{\ov{F}(x)} \ge \max(c^{(1)}, c^{(2)}).
\end{equation}
Now we could use (\ref{up1}) and obtain the upper bound
 \begin{equation}\label{up2}
\limsup_{x\to\infty} 
\frac{{\mathbf P} (Z>x)}{\ov{F}(x)} \le K_0,
\end{equation}
for some positive constant $K_0$. However, we need an explicit representation for
events leading to large values of $Z$. For that, we find it convenient to work with
a more ``rough'' upper bound than (\ref{up1}). Namely, take an arbitrary positive
number  $T>0$ and define random variables 
$\widehat{\sigma}_n$ by the equalities 
\begin{displaymath}
\widehat{\sigma}_n =
\Sigma_n
{\bf I} (\Sigma_n > T) + \widetilde{\sigma}_n {\bf I} (\Sigma_n \le T) 
\end{displaymath}
where
$$
\Sigma_n = \sum_{i=1}^L \left(
\sigma_{(n-1)L+i}^{(1)} + \sigma_{(n-1)L+i}^{(2)}
\right) \ge \widetilde{\sigma}_n
$$
and ${\bf I}$ is the indicator function.  
Clearly, $\widehat{\sigma}_n \ge \widetilde{\sigma}_n$ a.s. Further, 
due to Properties \ref{p1} and \ref{p2} from the Appendix, the common
distribution function 
of the random variables $\widehat{\sigma}_n$ belongs to the class $\SS{\gamma}$,
and $\Pr (\widehat{\sigma}_1 > x) \sim C \overline{F}(x)$  for some
positive $C$. Indeed, the distribution of $\Sigma_n$ 
belongs to the class $\SS{\gamma}$ by Property 2. 
Since, for any $x>T$,  
$\widehat{\sigma}_n >x$ if and only if 
$\Sigma_n>x$, we conclude that 
$\Pr (\widehat{\sigma}_n>x) \sim \Pr (\Sigma_n>x)$ as $x\to\infty$ and, therefore, 
the distribution 
$\widehat{\sigma}_n$ also belongs to the class $\SS{\gamma}$, due to Property 1..

We know that $\E \widetilde{\sigma} < \E \widetilde{\tau}_1 =
L\E \tau_1$ for all sufficiently large $L$. Also,  $\widetilde{\phi}(\gamma) = \E e^{\gamma 
\widetilde{\sigma}_1} <1/\E e^{-\gamma \widetilde{\tau}_1}$ 
for all 
large $L$, from \eqref{eq:11} and \eqref{star}. Choose such an $L$.
Further,
$
\E \widehat{\sigma} \to \E \widetilde{\sigma}$
and
$ \E e^{\gamma \widehat{\sigma}_1} \to \E e^{\gamma \widetilde{\sigma}_1}$
as $T\to\infty$. Therefore we can take a sufficiently large
$T$ for the inequalities 
\begin{equation}\label{hathat}
  \E \widehat{\sigma}_1 < L \E \tau_1 \quad \text{and}
\quad 
\E e^{\gamma \widehat{\sigma}_1} <1/\E e^{-\gamma \widetilde{\tau}_1} 
\end{equation}
to hold. Then the single-server queueing system with i.i.d. inter-arrival times
$\widetilde{\tau}_n$ (see Subsection 2.1 for the definition) 
and service times $\widehat{\sigma}_n$ is stable.
Let $\widehat{W}$ denote  
the stationary waiting time in this system, 
\begin{displaymath}
\widehat{W} = \max_{n\ge 0} 
\sum_{-n}^{-1}
(\widehat{\sigma}_n - \widetilde{\tau}_n) < \infty \quad\text{a.s.}
\end{displaymath}
Then $\widehat{W}$ coincides in distribution with the
supremum of a random walk with increments
$\widehat{\sigma}_n - \widetilde{\tau}_n$. By  
Property 4 from the Appendix and from \eqref{hathat}, 
$\Pr (\widehat{W}>x) \sim C\overline{F}(x)$ for some $C>0$ 
and, by Property 2, the distribution of the random variable
 $\widehat{W}$ belongs to the class 
${\cal S}_{\gamma}$. 
{From} the monotonicity properties, 
\begin{equation}\label{up10}
Z \le \widehat{Z} \equiv \widehat{W}+ \sigma_0^{(1)}+\sigma_0^{(2)}
\end{equation}
where the increments in the right-hand side are mutually independent, 
and the tail  distribution
of each of them is asymptotically equivalent to $\overline{F}(x)$,
up to a multiplicative non-negative constant (where at least one of these
constants is strictly positive). 
By Property 2 from the Appendix, the distribution of 
$\widehat{Z}$ also belongs to the class ${\cal S}_{\gamma}$ and
\begin{eqnarray*}
\Pr (\widehat{Z}>x) 
&=& 
\Pr \left(
\bigcup_{i=1}^2
\{
\widehat{Z}>x, {\sigma}_0^{(i)} > x- h_1(x)
\}
\bigcup
\{
\widehat{Z}>x, \widehat{W} > x- h_1(x)
\} \right) + o(\overline{F}(x)) \\
&=&
\sum_{i=1}^{2}\Pr (\widehat{Z}>x,{\sigma}_0^{(i)}> x- h_1(x))
+\Pr (\widehat{Z}>x, \widehat{W} > x-h_1(x))+o(\overline{F}(x))
\end{eqnarray*}
for any function $h_1$ satisfying \eqref{hx}. Note that if
 $h_1$ and $h$ are two such functions, then 
\begin{equation}\label{equi}
 \Pr (\widehat{Z}>x,{\sigma}_0^{(i)} > x- h_1(x))
=
 \Pr (\widehat{Z}>x,{\sigma}_0^{(i)} > x- h(x)) +o(\overline{F}(x)).
\end{equation}

Make use of the following simple relations. Let
$A$,$B$ and $C$ be three events. 
If
$$ 
\Pr (A) = \Pr (A \cap B) + v, 
$$
then 
\begin{equation}\label{ABC}
\Pr (A \cap C ) =
\Pr (A \cap C \cap B) + \widehat{v}
\end{equation}
where $0\le \widehat{v} \le v$. 
In particular, if
$C \subseteq A$, then the last equality
in transformed into 
$$ 
\Pr (C) = \Pr (C \cap B) + \widehat{v}.
$$
Applying this to the events
$A = \{ \widehat{Z}>x \} $ and $C = \{ Z>x\}$, 
we arrive at the equation 
\begin{equation}\label{gx}
\Pr ({Z}>x) 
= 
\sum_{i=1}^{2}\Pr ({Z}>x, {\sigma}_0^{(i)}\ge x- h_1(x))
+\Pr ({Z}>x, \widehat{W} > x-h_1(x))+o(\overline{F}(x)).
\end{equation}

By Property 4 from the Appendix, for any $\epsilon >0$, there exists 
a sufficiently large  $N$
such that
\begin{eqnarray*}
\Pr ( \widehat{W}>x-h_1(x) ) &=& 
\Pr
\left( \bigcup_{j=1}^N 
\{ \widehat{W}>x-h_1(x), \widehat{\sigma}_{-j}-\tau_{-j}>x-h_1(x)
-h_2(x-h_1(x))\}
\right) + g_1(x)\\
&=& 
\Pr
\left( \bigcup_{j=1}^N 
\{ \widehat{W}>x-h_1(x), \widehat{\sigma}_{-j}>x-h_1(x)
-h_2(x-h_1(x))\}
\right) + g_2(x)\\
&=& 
\sum_{j=1}^N \Pr 
(\widehat{W}>x-h_1(x), \widehat{\sigma}_{-j} >x-h_1(x)-h_2(x-h_1(x)))
+ g_3(x)
\end{eqnarray*}
where $h_2$ is any function satisfying (6) and
$0 \le g_i(x) \le \varepsilon \overline{F}(x) + o(\overline{F}(x))$
for $i=1,2,3$.
The latter inequality means that 
$\limsup_{x\to\infty} {g_i(x)}/{\overline{F}(x)} \le \varepsilon$.

Using again (\ref{ABC}), this time with $A = \{ \widehat{W}>x-h_1(x)\}$
and $C = \{Z>x\}$, we get
$$
\Pr (Z>x,\widehat{W}>x-h_1(x))
=
\sum_{j=1}^N \Pr (Z>x,\widehat{W}>x-h_1(x)), \widehat{\sigma}_{-j} > x-h_1(x)-
h_2(x-h_1(x))) + g_4(x)
$$
where
$0 \le g_4(x) \le \varepsilon \overline{F}(x) + o(\overline{F}(x))$.

By choosing an appropriate $h_2$, the right-hand side of the latter inequality
may be made simpler, by the exclusion the inequality $\{ \widehat{W}>x-h_1(x)\}$. 
Indeed, put
$$
P(x) = \Pr (Z>x,\widehat{W} \le x-h_1(x)), \widehat{\sigma}_{-j} > x-h_1(x)-
h_2(x-h_1(x)))
$$
and note that
\begin{eqnarray*}
P(x) 
&\le &
\Pr (\widehat{Z}>x,\widehat{W} \le x-h_1(x)), \widehat{\sigma}_{-j} > x-h_1(x)-
h_2(x-h_1(x))) \\
&\le &
\Pr (\sigma_0^{(1)} + \sigma_0^{(2)} > h_1(x),
\widehat{\sigma}_{-j} > x-h_1(x)-
h_2(x-h_1(x))) \\
&=& (1+o(1)) (c^{(1)}+c^{(2)}) \overline{F}(h_1(x)) \cdot C \overline{F}
(x-h_1(x)-
h_2(x-h_1(x)))\\
&=&
(1+o(1)) \widehat{C} e^{\gamma (h_1(x)+h_2(x-h_1(x)))}
\overline{F}(h_1(x))\overline{F}(x)
\end{eqnarray*}
where the last equality follows from the remark after formula (6) and
$\widehat{C} = C(c^{(1)}+c^{(2)})$. From
$\int_0^{\infty} e^{\gamma t} dF(t) < \infty$, we get $e^{\gamma t}\overline{F}(t)
\to 0$ as $t\to\infty$ and, therefore, $\overline{F}(h_1(x)) e^{\gamma h_1(x)} \to 0$
as $x\to\infty$. Thus, we may take $h_2$ so slowly increasing to infinity,
that
$\overline{F}(h_1(x)) e^{\gamma (h_1(x)+h_2(x-h_1(x)))} $ also tends to zero.
Then
$P(x) = o(\overline{F}(x))$ and, therefore,
$$
\Pr (Z>x,\widehat{W}>x-h_1(x))
=
\sum_{j=1}^N \Pr (Z>x, \widehat{\sigma}_{-j} > x-h_3(x)) + g_5(x)
$$
where
$0 \le g_5(x) \le \varepsilon \overline{F}(x) + o(\overline{F}(x))$
and the function $h_3(x) = h_1(x)+h_2(x-h_1(x))$ satisfies (6).

For all sufficiently large $x$, the events  $\{ \widehat{\sigma}_{-j} > x-h_3(x))\}$
 and $\{ \Sigma_{-j} > x-h_3(x))\}$ either occur or do not occur simultaneously.
Therefore $\Pr (Z>x, \widehat{\sigma}_{-j} > x-h_3(x)) =
\Pr (Z>x, \Sigma_{-j} > x-h_3(x))$.
Applying Property 2 from the Appendix to the random variables $\Sigma_{-j}$,
we get 
\begin{eqnarray*}
\Pr (Z>x, \widehat{\sigma}_{-j} > x-h_3(x)) 
&=&
\sum_{i=1}^2 \sum_{l=(j-1)L+1}^{jL} 
\Pr (Z>x,  \widehat{\sigma}_{-j} > x-h_3(x), \sigma_{-l}^{(i)}>
x - h_3(x) - h_4(x-h_3(x))) \\
&+& o(\overline{F}(x))
\end{eqnarray*}
where $h_4$ is any function satisfying (6).
Assuming that $h_4$ satisfies an extra condition (which is analogues
to that on $h_2$), we may exclude the inequality
$\{\widehat{\sigma}_{-j} > x-h_3(x)\}$ from the right-hand side
of the last relation.
Letting $h(x) = h_3(x)+h_4(x-h_3(x))$, we arrive at 
the equality 
$$
\Pr (Z>x, \widehat{\sigma}_{-j} > x-h_3(x)) =
\sum_{i=1}^2 \sum_{l=(j-1)L+1}^{jL} 
\Pr (Z>x,  \sigma_{-l}^{(i)}>
x - h(x)) + o(\overline{F}(x)).
$$

Substituting the relations obtained into formula (\ref{gx}) 
and using (\ref{equi}), we get finally:
\begin{equation}\label{g3}
 \Pr (Z>x) 
=
\sum_{i=1}^2 \sum_{j=0}^{NL}
\Pr (Z>x, \sigma_{-j}^{(i)} > x-h(x))
+ g(x)
\end{equation}
where
$0\le g(x) \le \varepsilon \overline{F}(x) + o(\overline{F}(x))$.

Now we study the asymptotics of each individual summand in the double sum in the right-hand side
of
(\ref{g3}). Let
$$
W_j^{(1)} = \max \left(
0, \sup_{n\ge 1} \sum_{i=1}^n (\sigma_{-j-i}^{(1)} - \tau_{-j-i})
\right)
$$
and notice that the distribution of
 $W_j^{(1)}$ does not depend on $j$. Further, let  
$$
Y_j^{(1)} = \max_{0\le k \le -j}
\left(
\sum_{i=-j+1}^{-j+k} \sigma_i^{(1)} + \sum_{i=-j+k}^0
\sigma_i^{(2)}
\right)
$$
and
$$
V_j^{(1)}
=
\max \left( 
\sup_{-j < n \le m < \infty}
\left( \sum_{-m}^{-n} \sigma_i^{(1)} + \sum_{-n}^0 \sigma_i^{(2)}
- \sum_{-m}^{-1} \tau_i \right),
\max_{0\le n \le m < -j} 
\left( \sum_{-m}^{-n} \sigma_i^{(1)} + \sum_{-n}^0 \sigma_i^{(2)}
- \sum_{-m}^{-1} \tau_i \right) \right).
$$
Then, for any $j\ge 0$,
$$
Z = \max \left( V_j^{(1)},
W_j^{(1)} + \sigma_{-j}^{(1)} + Y_j^{(1)}
- \sum_{i=-j}^{-1} \tau_i
\right)
$$
where the random variables $(W_j^{(1)},Y_j^{(1)},V_j^{(1)}, \sum_{-j}^{-1} \tau_i)$ 
mutually do not depend on $\sigma_{-j}^{(1)}$. 
Put $Q_j^{(1)}= W_j^{(1)}+Y_j^{(1)}-\sum_{-j}^{-1}
\tau_i$. 
For any $j=0, \ldots, NL$, 
\begin{eqnarray*}
\Pr (Z>x, \sigma_{-j}^{(1)}>x-h(x))
&= &
\Pr(Q_j^{(1)}+\sigma_{-j}^{(1)}>x,\sigma_{-j}^{(1)}>x-h(x))+o(\Pr(Z>x))\\
& = &
\int_0^{h(x)} \Pr (Q_j^{(1)} \in dt) \Pr (\sigma_{-j}^{(1)}>x-t)
+ o(\Pr(Z>x)) + o(\overline{F}(x))\\
&=& c^{(1)} \overline{F}(x) \int_0^{h(x)}
\Pr (Q_j^{(1)} \in dt) e^{-\gamma t} + o(\Pr(Z>x)+\overline{F}(x))\\
&=&
c^{(1)} \overline{F}(x) \E e^{\gamma W^{(1)}_0}
\E e^{\gamma Y_j^{(1)}} \left( \varphi_{\tau}(-\gamma) \right)^j
+ o(\Pr(Z>x)+\overline{F}(x))
\end{eqnarray*}
We clarify now each of four equalities above.
The first of them follows from
\begin{eqnarray*}
\Pr (Z>x, \sigma_{-j}^{(1)}>x-h(x))
&= &
\Pr(Q_j^{(1)}+\sigma_{-j}^{(1)}>x,\sigma_{-j}^{(1)}>x-h(x)) \\
&+& 
\Pr(V_j^{(1)}>x, Q_j^{(1)}+\sigma_{-j}^{(1)}\le x,\sigma_{-j}^{(1)}>x-h(x))
\end{eqnarray*}
where the second summand is not bigger than
\begin{eqnarray*}
\Pr(V_j^{(1)}>x, \sigma_{-j}^{(1)}>x-h(x)) 
&=&
\Pr(V_j^{(1)}>x)\Pr( \sigma_{-j}^{(1)}>x-h(x)) \\
&\le & \Pr(Z>x)\Pr( \sigma_{-j}^{(1)}>x-h(x)) = o(\Pr (Z>x)).
\end{eqnarray*}
Further, 
\begin{eqnarray*}
\Pr(Q_j^{(1)}+\sigma_{-j}^{(1)}>x,\sigma_{-j}^{(1)}>x-h(x)) 
&=&
\int_0^{h(x)} \Pr (Q_j^{(1)} \in dt) \Pr (\sigma_{-j}^{(1)}>x-t)\\ 
&+&
\Pr (Q_j^{(1)}>h(x)) \Pr (\sigma_{-j}^{(1)}>x-h(x))
\end{eqnarray*}
where
$$
\Pr (Q_j^{(1)}>h(x)) \le \Pr (Z>h(x))
$$
and
$$
\Pr (\sigma_{-j}^{(1)}>x-h(x)) \sim c^{(1)}\overline{F}(x-h(x))
\sim c^{(1)} e^{\gamma h(x)} \overline{F}(x).
$$
Since ${\mathbf E} e^{\gamma Z} < \infty$, we get 
$\Pr (Z>h(x)) e^{\gamma h(x)} \to 0$ when $x\to\infty$. Therefore, the second
equality also holds. The third equality follows from the uniform
equivalence 
(\ref{hx}) and from the assumptions of the theorem:
$$
\int_0^{h(x)} \Pr (Q_j^{(1)} \in dt) \Pr (\sigma_{-j}^{(1)}>x-t)
\sim
c^{(1)} \int_0^{h(x)} \Pr (Q_j^{(1)} \in dt) \overline{F}(x-t)
\sim c^{(1)}  \int_0^{h(x)} \Pr (Q_j^{(1)} \in dt) e^{\gamma t}\overline{F}(x).
$$
Finally, as $x\to\infty$,
$$
\int_0^{h(x)} \Pr (Q_j^{(1)} \in dt) e^{\gamma t}
 \rightarrow 
\int_0^{\infty} \Pr (Q_j^{(1)} \in dt) e^{\gamma t}
= {\mathbf E} e^{\gamma Q_j^{(1)}},
$$
and the last equality follows from the mutual independence of
the summands in 
$Q_j^{(1)}$.

Hence, 
\begin{equation}\label{odin}
\sum_{j=0}^{NL} 
\Pr (Z>x, \sigma_{-j}^{(1)}>x-h(x))
= (1+o(1))c^{(1)} \overline{F}(x) \E e^{\gamma W^{(1)}_0}
\sum_{j=0}^{NL} \E e^{\gamma Y_j^{(1)}} \left( \varphi_{\tau}(-\gamma) \right)^j
+ o(\Pr(Z>x)).
\end{equation}

Similarly, for any $j=0,1,2,\ldots$, the random variable $Z$ 
may be represented as
$$
Z= \max \left(V_J^{(2)}, 
Y_j^{(2)} + \sigma_{-j}^{(2)} + \sum_{i=-j+1}^0 \sigma_i^{(2)} -
\sum_{i=-j}^{-1} \tau_i \right)
$$
where 
$$
Y_j^{(2)}= \sup_{m\ge n \ge -j}
\left( \sum_{-m}^{-n} \sigma_i^{(1)} + \sum_{-n}^{-j-1}
\sigma_i^{(2)} - \sum_{-m}^{-j-1} \tau_i \right)
$$
(and the distribution of $Y_j^{(2)}$ does not depend on $j$),
$$
V_j^{(2)} = 
\sup_{m\ge -j} \max_{0\le n<-j}
\left( \sum_{-m}^{-n} \sigma_i^{(1)} + \sum_{-n}^{0}
\sigma_i^{(2)} - \sum_{-m}^{-1} \tau_i \right),
$$
and random variables $(Y_j^{(2)}, \sum_{i=-j+1}^0 \sigma_i^{(2)} -
\sum_{i=-j}^{-1} \tau_i, V_j^{(2)})$ are mutually independent of  
$\sigma_{-j}^{(2)}$. Then (with 
$Q_j^{(2)}=  Y_j^{(2)} + \sum_{i=-j+1}^0 \sigma_i^{(2)} -
\sum_{i=-j}^{-1} \tau_i $  )
\begin{eqnarray*}
\Pr (Z>x, \sigma_{-j}^{(2)}>x-h(x))
&= &
\Pr (Q_j^{(2)}+ \sigma_{-j}^{(2)}>x, \sigma_{-j}^{(2)}>x-h(x))
+ o(\Pr(Z>x))\\
&=&
\int_0^{h(x)} \Pr (Q_j^{(2)}\in dt) \Pr (\sigma_{-j}^{(2)}>x-t)
+o(\Pr(Z>x)+\overline{F}(x))\\
&= & c^{(2)} \overline{F}(x) 
\int_0^{h(x)}\Pr (Q_j^{(2)}\in dt) e^{\gamma t}
+o(\Pr(Z>x)+\overline{F}(x))\\
&= & c^{(2)} \overline{F}(x) 
\E e^{\gamma Y_0^{(2)}}\left(
\varphi_{(2)}(\gamma)\right)^{j-1}
\left(\varphi_{\tau}(-\gamma)\right)^{j-1}
+ o(\Pr(Z>x)+\overline{F}(x))
\end{eqnarray*}
(where the arguments are similar to those which were used to obtain 
the asymptotics for 
$\Pr (Z>x, \sigma_{-j}^{(2)}>x-h(x))$). 
Thus,
\begin{equation}\label{dva}
\sum_{j=0}^{NL} \Pr (Z>x, \sigma_{-j}^{(2)}>x-h(x))
= (1+o(1))
c^{(2)} \overline{F}(x) \E e^{\gamma Y_0^{(2)}}\frac{1-R^{NL}_2}{1-R_2}
+ o(\Pr(Z>x))
\end{equation}
where $R_2= \varphi_{(2)}(\gamma) \varphi_{\tau}(-\gamma) <1$.

Putting together (\ref{g3}), (\ref{odin}) and (\ref{dva}), we get:
\begin{eqnarray*}
 \Pr (Z>x) (1+o(1)) &=& (1+o(1))
\overline{F}(x)\\
&\times &\left(
c^{(1)} 
\E e^{\gamma W^{(1)}_0}
\sum_{j=0}^{NL} \E e^{\gamma Y_j^{(1)}} \left( \varphi_{\tau}(-\gamma) \right)^j
+
c^{(2)} \E e^{\gamma Y_0^{(2)}}\frac{1-R_2^{NL}}{1-R_2}\right) + g_3(x).
\end{eqnarray*}
Therefore,
$$
\limsup_{x\to\infty}
\frac{\Pr(Z>x)}{\overline{F}(x)}
\le 
c^{(1)} 
\E e^{\gamma W^{(1)}_0}
\sum_{j=0}^{\infty} \E e^{\gamma Y_j^{(1)}} \left( \varphi_{\tau}(-\gamma) \right)^j
+
c^{(2)} \E e^{\gamma Y_0^{(2)}}\frac{1}{1-R_2} + \varepsilon
$$
and
$$
\liminf_{x\to\infty}
\frac{\Pr(Z>x)}{\overline{F}(x)}
\ge 
c^{(1)} 
\E e^{\gamma W^{(1)}_0}
\sum_{j=0}^{NL} \E e^{\gamma Y_j^{(1)}} \left( \varphi_{\tau}(-\gamma) \right)^j
+
c^{(2)} \E e^{\gamma Y_0^{(2)}}\frac{1-R_2^{NL}}{1-R_2}
$$
for any positive
$\varepsilon$ (and, respectively,
for any positive integer $N$).
Letting $\varepsilon$ to zero, we obtain finally: 
\begin{equation}\label{zvezda}
 \Pr (Z>x) \sim
\overline{F}(x)
\left(
c^{(1)} 
\E e^{\gamma W^{(1)}_0}
\sum_{j=0}^{\infty} \E e^{\gamma Y_j^{(1)}} \left( \varphi_{\tau}(-\gamma) \right)^j
+
c^{(2)} \E e^{\gamma Y_0^{(2)}}\frac{1}{1-R_2}\right).
\end{equation}

We prove now the bounds (\ref{up-low}).
For this, we use repeatedly the following relations:
for any finite or countable collection of random variables $X_i$,  
$$
\sup_i {\mathbf E} e^{X_i} \le {\mathbf E} e^{\sup_i X_i} \le \sum_{i} 
{\mathbf E} e^{X_i}.
$$

First,
\begin{displaymath}
 1 \le  \E e^{\gamma W_0^{(1)}} 
\le  \sum_{n=0}^{\infty}
\left(
\E e^{\gamma (\sigma_1^{(1)}-\tau_1)}\right)^n 
=
\frac{1}{1-R_1}
\end{displaymath}
where again $R_1 =\varphi_{(1)}(\gamma) \varphi_{\tau}(-\gamma) <1$.
Further, with $\varphi (\gamma) = \max (\varphi_{(1)}(\gamma),
\varphi_{(2)}(\gamma))$,
\begin{eqnarray*}
{\mathbf E} e^{\gamma Y_j^{(1)}} 
&\ge &
\max_{0\le k \le -j}
{\mathbf E} \exp
\left(\gamma 
\sum_{i=-j+1}^{-j+k} \sigma_i^{(1)} + \gamma \sum_{i=-j+k}^0
\sigma_i^{(2)}
\right)\\
&=&
\max_{0\le k \le -j} \varphi_{(1)}^k (\gamma) \varphi_{(2)}^{j+1-k}(\gamma)
=\varphi^j(\gamma) \varphi_{(2)}(\gamma) 
\end{eqnarray*}
and
\begin{eqnarray*}
\E e^{\gamma Y_j^{(1)}} &\le &
\sum_{k=0}^{-j}
{\mathbf E} \exp
\left(\gamma 
\sum_{i=-j+1}^{-j+k} \sigma_i^{(1)} + \gamma \sum_{i=-j+k}^0
\sigma_i^{(2)}
\right)\\
&=&
\varphi_{(2)}(\gamma) \sum_{i=0}^{j}
\varphi_{(1)}^i(\gamma) \varphi_{(2)}^{j-i}(\gamma).
\end{eqnarray*}
Since $R=\varphi (\gamma ) \varphi_{\tau} (-\gamma)$, 
\begin{eqnarray*}
\frac{\varphi_{(2)}(\gamma)}{1-R} \le
\sum_{j=0}^{\infty} \E e^{\gamma Y_j^{(1)}} \varphi_{\tau}(-\gamma)^j \le 
\phi_{(2)}(\gamma) 
\cdot
\frac{1}{(1-R_1)(1-R_2)}.
\end{eqnarray*}
Similarly,
\begin{eqnarray*}
{\varphi_{(1)}(\gamma)} \le \E e^{\gamma Y^{(2)}_0} \le
\phi_{(1)}(\gamma) 
\cdot
\frac{1}{(1-R_1)(1-R_2)}.
\end{eqnarray*}
Substituting all obtained inequalities into (\ref{zvezda}),
we arrive at (\ref{up-low}).

\section{Appendix}

\subsection{Properties of distributions from the class $\SS{\beta}$, $\beta >0$.}

We present here a number of known properties (Properties 1--3)
of the class $\SS{\beta}$ with $\beta >0$ -- see, for example,
\cite{PAKES} and the comments in \cite{FZ}, and also Property 
4 which follows from results of \cite{FZ}. 

\begin{property} \label{p1} (Closure of the class
 $\SS{\beta}$ with respect to the tail equivalence) \\ 
If $F\in \SS{\beta}$
and $\overline{F}(x) \sim c \overline{G}(x)$ for some constant
$c \in (0,\infty )$, then $G\in \SS{\beta}$.
If, in particular,  random variables $X$ and $Y$ are independent, $Y$ is
a.s. non-negative and. $F(x) = {\bf P} (X\le x) \in \SS{\beta}$,
then
$$
\overline{G}(x) = {\bf P} (X-Y > x) \sim \int_0^{h(x)} 
G(dt) \overline{F}(x+t) \sim \overline{F}(x)
\int_0^{h(x)} G(dt) e^{-\beta t} \sim \overline{F}(x) {\mathbf E} e^{- \beta Y}
$$ 
as
$x\to\infty$ and, therefore,  
$G \in \SS{\beta}$.
\end{property}

The following more general result also holds.
\begin{property}\label{p2}
  Assume that $F\in\SS{\beta}$ for some  
  $\beta\ge0$.  Assume also that random variables  $X_i$, $i=1,\ldots,n$
are mutually independent and their distribution functions~$F_i$ satisfy the
relations   
  ${\mathbf P}(X_i>x)=\ov{F_i}(x)\sim{}c_i\ov{F}(x)$ 
as $x\to\infty$, for some  $c_i\ge0$,
$\sum_1^n c^{(i)}>0$.  Then
  $\phi_{i} = {\mathbf E} e^{\beta X_i}<\infty$ for all 
$i=1,\dots,n$. The distribution of the sum $\sum_{i=1}^nX_i$ 
also belongs to the class $ \SS{\beta}$ and
\begin{eqnarray*}
    {\bf P} (\sum_{i=1}^nX_i >x)
&\sim &
\sum_{j=1}^n {\bf P}
\left(\bigcup_{j=1}^n \{
\sum_{i\ne j} X_i \le h(x), \sum_{i=1}^nX_i >x
\}
\right)\\
&\sim &
\sum_{j=1}^n {\bf P}
\left(\sum_{i\ne j} X_i \le h(x), \sum_{i=1}^nX_i >x
\right)\\
&\sim &
\sum_{j=1}^n {\bf P}
\left(\bigcup_{j=1}^n \{
 X_j > x-h(x), \sum_{i=1}^nX_i >x
\}
\right)\\
&\sim &
\sum_{j=1}^n {\bf P}
\left(X_j > x-h(x), \sum_{i=1}^nX_i >x
\right)
    \sim
    \prod_{i=1}^n\phi_{i}\sum_{i=1}^n\frac{c_i}{\phi_{i}}
    \;\overline{F}(x)
  \end{eqnarray*}
where $h(x)$ is any function satisfying (\ref{hx}).
\end{property}

\begin{property} \label{p4}
 Assume that $F\in\SS{\beta}$ for some  
  $\beta\ge0$.  Assume that random variables 
$V,\xi$ and $\eta$ are mutually independent,
$\eta\ge 0$ a.s., $\Pr (V >x) \sim
c_1 \overline{F}(x)$ and $
\Pr (\xi>x) \sim c_2 \overline{F}(x)$ where $c_1\ge 0$ and $c_2>0$.
Then, for any function $h$ satisfying 
(\ref{hx}), 
\begin{eqnarray*}
\Pr (V+\xi-\eta>x, V\le h(x))
&\sim &
\Pr (V+\xi-\eta>x, V-\eta \le h(x))\\
&\sim &
\Pr (V+\xi-\eta>x, \xi-\eta \ge x-h(x))\\
&\sim &
\Pr (V+\xi-\eta>x, \xi \ge x-h(x))\\
&\sim &
c_2\overline{F}(x) \E e^{\beta V} \E e^{-\beta \eta}.
\end{eqnarray*}
\end{property}

\begin{property} \label{p5}
Consider a sequence of i.i.d. random variables
 $\{ X_n\}$
with a common distribution function $F$ and assume that 
$\E X_i = -a < 0$, $F \in {\cal S}_{\beta}$ and
$\E e^{\beta X_1}<1.$ Let $M_k = \max_{0\le n \le k} \sum_{i=1}^n X_i$ 
and 
$M = \sup_{n\ge 0} \sum_{i=1}^n X_i.$ It follows from [5,  1 and Remark 3]
that
$$
\lim_{x\to\infty} \frac{\Pr (M>x)}{\overline{F}(x)}
\ = \ \frac{{\mathbf E} e^{\beta M}}{1-{\mathbf E} e^{\beta X_1}}
$$
(see also \cite{BD})
and, moreover, for any (as small as possible)
$\epsilon \in (0,1)$, there exists a sufficiently large
$N$ such that, as $x\to\infty$ and for any function
$h(x)$ satisfying \eqref{hx},
\begin{eqnarray*} 
\Pr (M>x) &\ge & 
\Pr 
\left( 
M>x, \bigcup_{n=1}^N \{ X_n>x-h(x)\}
\right) + o(\overline{F}(x))\\
&=&
\sum_{n=1}^N \Pr (M>x, X_n>x-h(x)) + o(\overline{F}(x))\\
&\ge &
\Pr (M>x) + o(\overline{F}(x)) -\varepsilon \overline{F}(x)
\end{eqnarray*}
(recall that the notation $f(x) \ge g(x) + o(f(x))$
means $\limsup_{x\to\infty} g(x)/f(x) \le 1$; in our case 
$o(\overline{F}(x))=o(\Pr(M>x))$). 

If we assume in addition that the random variables 
$X_n$ are represented as differences $X_n=\xi_n-\eta_n$ where
$\{\xi_n\}$ and $\{\eta_n\}$ are two mutually independent sequences of
i.i.d. random variables and random variables $\xi_n$ are non-negative a.s.,
then the relations above stay valid if we replace the events
$\{X_n>x-h(x)\}$ by $\{\xi_n > x-h(x)\}.$ 
\end{property}

\subsection{A simple proof of Theorem 1 by the use of upper and lower bounds from
Subsection 2.1.}
\label{trivial}

Recall again that the stationary waiting time in a single-server queue
with service times $\sigma_n$ and inter-arrival times $\tau_n$
has the same distribution as the supremum $M= \sup_n S_n$ 
of a random walk
$S_n=\sum_1^n X_i$ with increments $X_n = \sigma_n-\tau_n$.
Apply the following well-known result for a random walk 
$S_n = \sum_{i=1}^n X_i$ with negative drift:
\begin{theorem}\label{basic}
If
\begin{displaymath}
\lambda_0 = \sup \{ \lambda \ : \ \phi_{X_1} (\lambda ) \le 1 \} >0,
\end{displaymath}
then
\begin{displaymath}
-\ln {\mathbf P} (M>x) \sim \lambda_0 x.
\end{displaymath}
\end{theorem}
\begin{rem} 
We are aware of only one publication (\cite{GOW}, p. 17) where Theorem
\ref{basic} is formulated without extra assumptions. Usually authors
assume in addition (see, for example, 
[10, Section 21, Theorem 11]) the so-called Cramer condition
\begin{equation}\label{extra}
\phi_{X_1}(\lambda_0)=1 \quad\text{and} \quad
d={\mathbf E} X_1 e^{\lambda_0 X_1} < \infty 
\end{equation} 
or even stronger conditions. The theorem may be obtained also
as a corollary of more general results, for instance, from
 \cite{BM}. We provide in Subsection 3.3 a methodological
comment on how one can prove the general result of Theorem \ref{basic}
given that it is already proven under the conditions \eqref{extra}.  
\end{rem}
It follows from Theorem \ref{basic} that 
\begin{displaymath}
-\ln {\mathbf P} (W^{(i)}>x) \sim \gamma^{(i)} x,
\end{displaymath}
and since $Z^{(i)}\ge W^{(i)}$ a.s.,
\begin{equation} \label{oneone}
\limsup_{x\to\infty} \frac{-\ln {\mathbf P}(Z>x)}{x} \le \gamma.
\end{equation}

On the other hand, for any
$\lambda > 0$,
\begin{eqnarray*}
\widetilde{\phi}(\lambda) &\equiv& 
{\mathbf E} e^{\lambda\widetilde{\sigma}_1} 
\le \sum_{j=1}^L {\mathbf E} 
\exp \left(
\sum_{i=1}^j \sigma_i^{(1)}
+ \sum_{i=j}^L \sigma_i^{(2)}
\right)\\
&=& \sum_{j=1}^L \phi_{(1)}^j(\lambda)  
\phi_{(2)}^{L-j+1}(\lambda)
\equiv \phi_{*}(\lambda).
\end{eqnarray*}
Using again the notation $\phi (\lambda) = \max (\phi_{(1)}(\lambda),
\phi_{(2)}(\lambda))$, we get 
\begin{displaymath}
\min ((\phi_{(1)}(\lambda),
\phi_{(2)}(\lambda)) \cdot \phi^L(\lambda)
\le \widetilde{\phi}(\lambda)  
\le \phi_{*}(\lambda)
\le (L+1) \phi^{L+1}(\lambda)
\end{displaymath}
and therefore 
\begin{equation}\label{star}
\left(
\widetilde{\phi}(\lambda)
\right)^{1/L} \rightarrow \phi (\lambda) \quad
\text{and} \quad 
\left(
\phi_{*}(\lambda)
\right)^{1/L} \rightarrow \phi (\lambda)
 \quad\text{as}
\quad L\to\infty.
\end{equation}
Let $\widetilde{\gamma} = \sup \{ \lambda \ : \
\widetilde{\phi}(\lambda) \phi_{\tau}^L(-\lambda) \le 1 \}$.
Since $\widetilde{\sigma}_1 \ge \max \left(\sum_1^L \sigma^{(1)}_i,
\sum_1^L \sigma^{(2)}_i\right)$, we may conclude that
$\gamma \ge \widetilde{\gamma}$. 
>From \eqref{star},
$\widetilde{\gamma}\to\gamma$.
By Theorem \ref{basic}, for any sufficiently large $L$,
$$
-\ln {\mathbf P} (\widetilde{W}>x) \sim \widetilde{\gamma}x.
$$
>From ${\mathbf E} \exp
\left(
{\gamma \widetilde{\sigma}_0}
\right)<\infty$,
we get 
$$
-\ln {\mathbf P} (\widetilde{Z}>x) \sim \widetilde{\gamma}x.
$$
Letting $L$ to infinity, we obtain
\begin{equation} \label{twotwo}
\liminf_{x\to\infty} \frac{-\ln {\mathbf P}(Z>x)}{x} \ge \gamma.
\end{equation}
The statement of Theorem 1 follows now from the inequalities (\ref{oneone})
and (\ref{twotwo}).

\begin{rem}
A natural analogue of Theorem 1 
holds for a tandem of any finite number of queues, with a similar
proof.
\end{rem}

\subsection{A comment on a proof of Theorem \ref{basic}.}

Assume that the statement of the theorem has been already proved
under the additional assumptions \eqref{extra}.
Note that there are several versions of such a proof (by 
the use of, say, (a) martingale techniques, (b) exponential
change of measure and elements of the renewal theory, etc.)

Assume now that conditions \eqref{extra} do not hold.
For any $r>0$, define random variables 
\begin{displaymath}
X_{n,+r} = \max (X_n, -r) 
\quad\text{and}\quad X_{n,-r} = \min (X_n, r).
\end{displaymath}
Denote the corresponding sums, maxima and moments by
$S_{n,+r},S_{n,-r}, M_{+r}, M_{-r}$, $\phi_{X,+r}$ and $\phi_{X,-r}$ 
where $M_{+r} \ge M \ge M_{-r}$ a.s.
For all sufficiently large values of $r$, the maximum  $M_{+r}$
is a.s. finite. From the 
monotonicity and continuity of $\phi_{X,+r}$ and $\phi_{X,-r}$ as
functions of $r$ and from the boundedness from above of random variables
$X_{n,-r}$, it follows, firstly, that the roots $\lambda_{+r} < \lambda_0 < \lambda_{-r}$
of equations $\phi_{X,+r}(\lambda_{+r})=1$ and
$\phi_{X,-r}(\lambda_{-r})=1$ exist for any $r$, the corresponding derivatives
are finite and therefore
\begin{displaymath}
-\ln {\mathbf P} (M_{+r}>x) \sim \lambda_{+r}x 
\quad\text{and}\quad 
-\ln {\mathbf P} (M_{-r}>x) \sim \lambda_{-r}x, 
\end{displaymath}
and, second, both $\lambda_{+r}$ and $\lambda_{-r}$
converge to $\lambda_0$ as $r\to\infty$. 

If $\phi_{X}(\lambda_0)<1$, then there is  $r<\infty$
such that $\phi_{X,+r} (\lambda_0) =1$ and $\phi_{X,+r}^{'}
(\lambda_0) < \infty$. Then, for this $r$,
$$
-\ln {\bf P} (M_{+r} >x) \sim \lambda_0x,
$$
and the statement of Theorem \ref{basic} follows.

The case $\phi_{X}(\lambda_0)=1$ and $d=\infty$ is left for a reader.

The author would like to thank Stan Zachary for improving the style of the
English translation, and the referee for a number of important comments
and remarks.

{\it Remark} (added at the proofreading): in the paper \cite{Lel}, the author develops 
the approach
from \cite{BF95} to obtain 
the logarithmic asymptotics in a wide class of stochastic networks.


\newpage

\begin{thebibliography}{20}

\bibitem{Lo}
Loynes, R. M. 
``The stability of a system of queues in series,'' 
{\em Math. Proc. Cambridge Phil.
Soc.}, 1964, vol. 60, no.3, pp. 569--574.

\bibitem{BF94}
Baccelli, F. and Foss, S.
``Ergodicity of Jackson-Type Queueing Networks'',  
{\em Queueing Systems}, 1994, vol. 17, no. 1, pp. 5--72.

\bibitem{BPT}
Bertsimas, D., Paschalidis, I., and Tsitsiklis, J. 
``On the large deviation behaviour in acyclic networks of
G/G/1 queues'',
{\em Ann. \ Appl. \ Prob.}, 1998, vol. 8, no. 4, pp. 1027--1069. 

\bibitem{GAN} 
Ganesh, A.\ J.
``Large deviations of the sojourn time for queues in series'',
{\em Ann. \ Oper. \ Res.}, 1998, vol. 79, no. 1, pp. 3--26.


\bibitem{FZ} Zachary, S. and Foss, S.\ G.  
``On the exact distributional asymptotics for the suprumem of a 
random walk with increments
in a class of light-tailed distributions'', 
{\em Siberian Math. J.}, 2006, 
v. 47 , np. 6, pp. 1034--1041.


\bibitem{BF04}  Baccelli, F. and Foss, S. 
``Moments and tails in monotone-separable stochastic networks'', 
{\em Ann.\ Appl.\ Prob}, 2004, vol. 14, no. 3,  pp. 612--650.

\bibitem{BF95} Baccelli, F. and Foss, S. 
``On the Saturation Rule for the Stability of Queues''. 
{\em J.\ Appl.\ Prob.}, 1995, vol. 32, no. 2,  pp. 494--507

\bibitem{GOW}
Ganesh, A.\ J., O Connell, N., Wischik, D.
``{\em Big queues}'', Lecture Notes in Mathematics, vol. 1838,
Springer, 2004.


\bibitem{BOR} Borovkov, A.\ A. 
``{\em Stochastic Processes in Queueing Theory}''.
Wiley, New York, 1976.

\bibitem{BM} Borovkov, A.\ A.\ and Mogul'skii A.\ A.
``The second deviations function and the asymptotic problems of 
renewal and hitting the boundary for multidimensional random walks''
{\em Siberian Math. J.}, 1996, 
v. 37, no. 4, pp. 647--682.


\bibitem{PAKES}
  Pakes, A.
  ``On the tails of waiting time distributions.''
  {\em J.\ Appl.\ Prob.}, 1975, vol. 7, no.4, pp. 745--789.


\bibitem{BD} 
 Bertoin, J.\ and Doney, R.\ A. 
  ``Some asymptotic results for transient random walks'', 
  {\em Adv.\ Appl.\ Prob.}, 1996, vol. 28, no.1, pp. 207--226.


\bibitem{Lel} Lelarge, M. 
``Tail asymptotics for monotone-separable networks'',
{\em J. Appl. Prob.}, 2007, vol. 44, no. 2, pp. 306--320.

\end{thebibliography}
\end{document}